\documentclass{article}%
\usepackage{amsfonts}
\usepackage{graphicx}
\usepackage{amsmath}
\usepackage{amssymb}%
\providecommand{\U}[1]{\protect\rule{.1in}{.1in}}
\newtheorem{theorem}{Theorem}

\newtheorem{corollary}[theorem]{Corollary}

\newtheorem{lemma}[theorem]{Lemma}

\newtheorem{proposition}[theorem]{Proposition}
\newtheorem{remark}[theorem]{Remark}

\begin{document}

\begin{center}
{\LARGE Almost Sure Invariance Principles\ via Martingale Approximation}

\bigskip

\bigskip Florence Merlev\`{e}de$^{a}$, Costel Peligrad$^{b}$ and Magda
Peligrad\footnote{Supported in part by a Charles Phelps Taft Memorial Fund
grant and a NSA\ grant.}$^{c}$
\end{center}

$^{a}$ Universit\'{e} Paris Est, Laboratoire de math\'{e}matiques, UMR 8050
CNRS, B\^{a}timent Copernic, 5 Boulevard Descartes, 77435 Champs-Sur-Marne,
FRANCE. E-mail: florence.merlevede@univ-mlv.fr

$^{b}$ Department of Mathematical Sciences, University of Cincinnati, PO Box
210025, Cincinnati, Oh 45221-0025, USA. E-mail address: peligrc@ucmail.uc.edu

$^{c}$ Department of Mathematical Sciences, University of Cincinnati, PO Box
210025, Cincinnati, Oh 45221-0025, USA. E-mail address: peligrm@ucmail.uc.edu

\begin{center}
\bigskip

Abstract
\end{center}

In this paper we estimate the rest of the approximation of a stationary
process by a martingale in terms of the projections of partial sums. Then,
based on this estimate, we obtain almost sure approximation of partial sums by
a martingale with stationary differences. The results are exploited to further
investigate the central limit theorem and its invariance principle started at
a point, as well as the law of the iterated logarithm via almost sure
approximation with a Brownian motion, improving the results available in the
literature. The conditions are well suited for a variety of examples; they are
easy to verify, for instance, for linear processes and functions of Bernoulli shifts.

Key words: martingale approximation, quenched CLT, normal Markov chains,
functional CLT, law of the iterated logarithm, almost sure approximation.

\bigskip

AMS 2000 Subject Classification: Primary 60F05, 60F15, 60J05.

\section{I\textbf{ntroduction and notations}}

In recent years there has been an intense effort towards a better
understanding of the structure and asymptotic behavior of stochastic
processes. For processes with short memory there are two basic techniques:
approximation with independent random variables or with martingales. Each of
these methods have its own strength. On one hand the classes that can be
treated by coupling with an independent sequence exhibit faster rates of
convergence in various limit theorems; on the other hand the class of
processes that can be treated by a martingale approximation is larger. There
are plenty of processes that benefit from approximation with a martingale.
Examples are: linear processes with martingale innovations, functions of
linear processes, reversible Markov chains, normal Markov chains, various
dynamic systems, discrete Fourier transform of general stationary sequences. A
martingale approximation provides important information about these
structures, and since martingales can be embedded into Brownian motion, they
satisfy the functional central limit theorem started at a point, and the law
of the iterated logarithm. Moreover, martingale approximation provides a
simple and unified approach to asymptotic results for many dependence
structures. For all these reasons, in recent years martingale approximation,
"coupling with a martingale", has gained a prominent role in analyzing
dependent data. This is also due to important developments by Liverani (1996),
Maxwell and Woodroofe (2000), Derriennic and Lin (2001 a), Wu and Woodroofe
(2004) and recent developments by Peligrad and Utev (2005), Peligrad, Utev and
Wu (2007), Merlev\`{e}de and Peligrad (2006), Peligrad and Wu (2010) among
others. Many of these new results, originally designed for Markov operators,
(see Kipnis and Varadhan, 1986; Derriennic and Lin, 2007, for a survey) have
made their way into limit theorems for stochastic processes.

So far this method has been shown to be well suited to transport from the
martingale to the stationary process either the conditional central limit
theorem or conditional invariance principle in mean. As a matter of fact,
papers by Dedecker-Merlev\`{e}de-Voln\'{y} (2006), Voln\'{y} (2007), Zhao and
Woodroofe (2007 a), Gordin and Peligrad (2010), obtain characterizations of
stochastic processes that can be approximated by martingales in quadratic
mean. These results are useful for treating evolutions in "annealed" media.

In this paper we address the question of almost sure approximation of partial
sums by a martingale. These results are useful for obtaining almost sure limit
theorems for dependent sequences and also limit theorems started at a point.
Limit theorems for stochastic processes that do not start from equilibrium is
timely and motivated by evolutions in quenched random environment. Moreover
recent discoveries by Voln\'{y} and Woodroofe (2010) show that many of the
central limit theorems satisfied by classes of stochastic processes in
equilibrium, fail to hold when the processes are started from a point, so, new
sharp sufficient conditions should be pointed out for the validity of these
types of results. Recent steps in this direction are papers by Zhao and
Woodroofe (2007 b), Cuny (2009 a and b), Cuny and Peligrad (2009).

The technical challenge is to estimate the rest of approximation of partial
sums by a martingale which leads to almost sure results, ranging from the
almost sure central limit theorems started at a point, almost sure
approximation with a Brownian motion and the law of the iterated logarithm.

We shall develop our results in the framework of stationary processes that can
be introduced in several equivalent ways.

We assume that $(\xi_{n})_{n\in\mathbb{Z}}$ denotes a stationary Markov chain
defined on a probability space $(\Omega,\mathcal{F},\mathbb{P})$ with values
in a measurable space. The marginal distribution and the transition kernel are
denoted by $\pi(A)=\mathbb{P}(\xi_{0}\in A)$ and $Q(\xi_{0},A)=\mathbb{P}%
(\xi_{1}\in A|\xi_{0})$. {Next let $\mathbb{L}_{0}^{2}(\pi)$ be the set of
functions on $S$ such that $\int f^{2}d\pi<\infty$ and $\int fd\pi=0,$ and for
a }${f}\in${$\mathbb{L}_{0}^{2}(\pi)$ denote $X_{i}=f(\xi_{i})$ }, $S_{n}%
=\sum\limits_{i=0}^{n-1}X_{i}$ (i.e. $S_{1}=X_{0}$, $S_{2}=X_{0}+X_{1},$
...).{ }In addition $Q$ denotes the operator on {$\mathbb{L}_{2}(\pi)$}
{acting via $(Qf)(x)=\int_{S}f(s)Q(x,ds)$}$.$ Denote by $\mathcal{F}_{k}$ the
$\sigma$--field generated by $\xi_{i}$ with $i\leq k$. For any integrable
variable $X$ we denote $\mathbb{E}_{k}(X)=\mathbb{E}(X|\mathcal{F}_{k}).$ In
our notation $\mathbb{E}_{0}(X_{1})=Qf(\xi_{0})=\mathbb{E}(X_{1}|\xi_{0}).$

Notice that any stationary sequence $(X_{k})_{k\in\mathbb{Z}}$ can be viewed
as a function of a Markov process $\xi_{k}=(X_{i};i\leq k),$ for the function
$g(\xi_{k})=X_{k}$.

The stationary stochastic processes may also be introduced in the following
alternative way. Let $T:\Omega\mapsto\Omega$ be a bijective bi-measurable
transformation preserving the probability. Let $\mathcal{F}_{0}$ be a $\sigma
$-algebra of $\mathcal{F}$ satisfying $\mathcal{F}_{0}\subseteq T^{-1}%
(\mathcal{F}_{0})$. We then define the nondecreasing filtration $(\mathcal{F}%
_{i})_{i\in\mathbb{Z}}$ by $\mathcal{F}_{i}=T^{-i}(\mathcal{F}_{0})$. Let
$X_{0}$ be a random variable which is $\mathcal{F}_{0}$-measurable. We define
the stationary sequence $(X_{i})_{i\in\mathbb{Z}}$ by $X_{i}=X_{0}\circ T^{i}$.

In this paper we shall use both frameworks. The variable $X_{0}$ will be
assumed centered at its mean, i.e. $\mathbb{E}(X_{0})=0$ and square integrable
$\mathbb{E}(X_{0}^{2})<\infty.$

The following notations will be frequently used. {We denote by }${||X}${$||$}
the norm in {$\mathbb{L}^{2}$}$(\Omega,\mathcal{F},\mathbb{P})$, the space of
square integrable functions. We shall also denote { by }${||X}${$||_{p}$} the
norm in {$\mathbb{L}^{p}$}$(\Omega,\mathcal{F},\mathbb{P}).$ For any two
positive sequences $a_{n}\ll b_{n},$ means that for a certain numerical
constant $C$, we have $a_{n}\leq Cb_{n}$ for all $n$; $[x]$ denotes the
largest integer smaller or equal to $x$. For the law of the iterated logarithm
we use the notation $\log_{2}n=\log(\log(\max(e,n)))$. The notation $a.s.$
means almost surely, while $\Rightarrow$ denotes convergence in distribution.

The main question addressed is to find sufficient projective conditions such
that there is a martingale $M_{n}$ with stationary differences such that
either%
\[
S_{n}-M_{n}=o(n^{1/2})\text{ a.s.}\,,
\]
or%
\[
S_{n}-M_{n}=o(n\log_{2}n)^{1/2}\text{ a.s.}\,
\]

These types of approximations are important to study for instance the limit
theorems stated at a point (quenched) and the law of the iterated logarithm
via almost sure approximation with a Brownian motion.

The "so called" quenched CLT, states that for any function $f$ continuous and bounded%

\begin{equation}
\mathbb{\ E}_{0}(f(S_{n}/\sqrt{n}))\rightarrow\mathbb{\ E}(f(cN))\text{
\ a.s.} \, , \label{CLT}%
\end{equation}
where $N$ is a standard normal variable and $c\ $\ is a certain positive
constant. By the quenched invariance principle we understand that for any
function $f$ continuous and bounded on $D[0,1]$ endowed with uniform topology
we have%

\begin{equation}
\mathbb{\ E}_{0}(f(S_{[nt]}/\sqrt{n}))\rightarrow\mathbb{\ E}(f(cW))\text{
\ a.s.} \label{FCLT}%
\end{equation}
where $W$ is the standard Brownian measure. We shall also refer to these types
of convergence also as almost sure convergence in distribution under
$\mathbb{P}_{0}$ a.s., where $\mathbb{P}_{0}(A)=\mathbb{P}(A|\mathcal{F}%
_{0}).$

This conditional form of the CLT is a stable type of convergence that makes
possible the change of measure with a majorating measure, as discussed in
Billingsley (1999), Rootz\'{e}n (1976), and Hall and Heyde (1980).

In the Markov chain setting the almost sure convergence in (\ref{CLT}) or
(\ref{FCLT})\ are presented in a slightly different terminology. Denote by
$\mathbb{P}^{x}$ and $\mathbb{E}^{x}$ the regular probability and conditional
expectation given $X_{0}=x$. In this context the quenched CLT\ is known under
the name of CLT started at a point i.e. the CLT or its functional form holds
for $\pi-$almost all $x\in S,$ under the measure $\mathbb{P}^{x}$.

Here is a short history of the quenched CLT under projective
criteria.\textbf{\ }A result in Borodin and Ibragimov (1994, Ch 4) states that
if $||\mathbb{E}_{0}(S_{n})||$ is bounded, then the CLT in its functional form
started at a point holds. Later work by Derriennic and Lin (2001 b) improved
on this result imposing the condition $||\mathbb{E}_{0}(S_{n})||\ll
n^{1/2-\epsilon}$ with $\epsilon>0$ (see also Rassoul-Agha and
Sepp\"{a}l\"{a}inen, 2008 and 2009). This condition was improved in Zhao and
Woodroofe (2008-a) and further improved by Cuny (2009 a) who imposed the
condition $||\mathbb{E}_{0}(S_{n})||\ll n^{1/2}(\log n)^{-2}(\log\log
n)^{-1-\delta}$ with $\delta>0$. A result in Cuny and Peligrad (2009) shows
that the condition $\sum_{k=1}^{\infty}||\mathbb{E}_{0}(X_{k})||/k^{1/2}%
<\infty$, is sufficient for (\ref{CLT}).

We shall prove here that the condition imposed to $||\mathbb{E}_{0}(S_{n})||$
can be improved, by requiring less restrictive conditions on the regularity of
$||\mathbb{E}_{0}(S_{n})||$ than the result in Cuny (2009 a)$.$ Then we shall
point out that the condition can be further weaken if we are interested in a
result for averages or if finite moments of order larger than $2$ are available.

To prove the law of the iterated logarithm we shall develop sufficient
conditions for almost sure approximation with a Brownian motion; that is we
shall redefine $X_{n},$ without changing its distribution, on a richer
probability space on which there exists a standard Brownian motion
$(W(t),t\geq0)$ such that for a certain positive constant $c>0,$\textit{ }%
\begin{equation}
S_{n}-W(cn)=o(n\log_{2}n)^{1/2}\text{ a.s.} \label{LIL}%
\end{equation}
Our method of proof is based on martingale approximation that is valid under
the Maxwell-Woodroofe condition:%
\begin{equation}
\Delta(X_{0})=\sum_{k=1}^{\infty}\frac{||\mathbb{E}_{0}(S_{k})||}{k^{3/2}%
}<\infty\text{ .} \label{MW}%
\end{equation}
The key tool in obtaining our results is the estimate of the rest of the
martingale approximation in terms of $||\mathbb{E}_{0}(S_{k})||$. We shall
establish in Section 2 that there is a unique martingale with stationary and
square integrable differences such that%

\begin{equation}
\frac{||S_{n}-M_{n}||}{n^{1/2}}\ll\sum_{k\geq n}\frac{||\mathbb{E}_{0}%
(S_{k})||}{k^{3/2}}\text{ .} \label{Approx}%
\end{equation}

We then further exploit the estimate (\ref{Approx}) to derive almost sure
martingale approximations of the types (\ref{AS1}) and (\ref{AS}) in Section 3.

Our paper is organized as follows: In Section 2 we present the martingale
approximation and estimate its rest. In Section 3 we present the almost sure
martingale approximation results. Section 4 is dedicated to almost sure
limiting results for the stationary processes. Section 5 points out some
applications. Several results involving maximal inequalities and several
technical lemmas are presented in the Appendix.

\section{Martingale approximation with rest}

\begin{proposition}
\label{approx}For any stationary sequence $(X_{k})_{k\in\mathbb{Z}}$ and
filtration $(\mathcal{F}_{k})_{k\in\mathbb{Z}}$ described above with
$\Delta(X_{0})<\infty$, there is a martingale $(M_{k})_{k\geq1}$ with
stationary and square integrable differences $(D_{k})_{k\in\mathbb{Z}}$
adapted to $(\mathcal{F}_{k+1})_{k\in\mathbb{Z}}$, $M_{n}=\sum_{i=0}%
^{n-1}D_{i}$, satisfying (\ref{Approx}).
\end{proposition}

To prove this proposition we need two preparatory lemmas. It is convenient to
use the notation
\begin{equation}
Y_{k}^{m}=\frac{1}{m}\mathbb{E}_{k}(X_{k+1}+...+X_{k+m})\text{ .} \label{defY}%
\end{equation}

We shall also use the following semi-norm notation. For a stationary process
$(X_{k})_{k\in\mathbb{Z}}$ define the semi-norm
\begin{equation}
||X_{0}||_{+}^{2}=\lim\sup_{n\rightarrow\infty}\frac{1}{n}\mathbb{E}(S_{n}%
^{2})\text{ .} \label{norm+}%
\end{equation}

\begin{lemma}
\label{estplus}Assume $||Y_{0}^{m}||_{+}\rightarrow0.$ Then, there is a
martingale $(M_{k})_{k\in\mathbb{Z}}$ with stationary and square integrable
differences adapted to $(\mathcal{F}_{k+1})_{k\in\mathbb{Z}}$ satisfying%
\[
\frac{||S_{n}-M_{n}||}{n^{1/2}}\ll\max_{1\leq l\leq n}\frac{||\mathbb{E}%
(S_{l}|\mathcal{F}_{0})||}{n^{1/2}}+||Y_{0}^{n}||_{+}\text{ .}%
\]

\end{lemma}

\textbf{Proof.} The construction of the martingale decomposition is based on
averages. It was introduced by Wu and Woodroofe (2004; see their definition 6
on the page 1677) and further developed in Zhao and Woodroofe (2008 b),
extending the construction in Heyde (1974) and Gordin and Lifshitz (1981); see
also Theorem 8.1 in Borodin and Ibragimov (1994), and Kipnis and Varadhan
(1986). We give the martingale construction with the estimation of the rest.

We introduce a parameter $m\geq1$ (kept fixed for the moment), and define the
stationary sequence of random variables:%
\[
\theta_{0}^{m}=\frac{1}{m}\sum_{i=1}^{m}\mathbb{E}_{0}(S_{i}),\text{ }%
\theta_{k}^{m}=\theta_{0}^{m}\circ T^{k}\text{ .}%
\]
Set
\begin{equation}
D_{k}^{m}=\theta_{k+1}^{m}-\mathbb{E}_{k}(\theta_{k+1}^{m})\text{ ; }M_{n}%
^{m}=\sum_{k=0}^{n-1}D_{k}^{m}\text{ .}\label{martd}%
\end{equation}
Then, $(D_{k}^{m})_{k\in\mathbb{Z}}$ is a stationary martingale difference,
$D_{k}^{m}$ is $\mathcal{F}_{k+1}$-measurable and sequence and $(M_{n}%
^{m})_{n\geq0}$ is a martingale. So we have%
\[
X_{k}=D_{k}^{m}+\theta_{k}^{m}-\theta_{k+1}^{m}+\frac{1}{m}\mathbb{E}%
_{k}(S_{k+m+1}-S_{k+1})\text{ }%
\]
and therefore%
\begin{align}
S_{k} &  =M_{k}^{m}+\theta_{0}^{m}-\theta_{k}^{m}+\sum\nolimits_{j=1}^{k}%
\frac{1}{m}\mathbb{E}_{j-1}(S_{j+m}-S_{j})\label{martdec}\\
&  =M_{k}^{m}+\theta_{0}^{m}-\theta_{k}^{m}+\overline{R}_{k}^{\,m}\text{
,}\nonumber
\end{align}
where we implemented the notation%
\[
\overline{R}_{k}^{\,m}=\sum\nolimits_{j=1}^{k}\frac{1}{m}\mathbb{E}%
_{j-1}(S_{j+m}-S_{j})\text{ .}%
\]
Observe that
\begin{equation}
\overline{R}_{k}^{\,m}=\sum_{j=0}^{k-1}Y_{j}^{m}.\label{rests}%
\end{equation}
With the notation
\begin{equation}
R_{k}^{m}=\theta_{0}^{m}-\theta_{k}^{m}+\overline{R}_{k}^{\,m}\,,\label{rest}%
\end{equation}
we have
\begin{equation}
S_{k}=M_{k}^{m}+R_{k}^{m}\text{ .}\label{deco}%
\end{equation}
Notice that
\begin{equation}
||S_{n}-M_{n}^{n}||\leq3\max_{1\leq i\leq n}||\mathbb{E}(S_{i}|\mathcal{F}%
_{0})||\text{ .}\label{estRn}%
\end{equation}
It was shown in Gordin and Peligrad (2009) that if $||Y_{0}^{m}||_{+}%
\rightarrow0,$ then $D_{0}^{n}$ converges in $\mathbb{L}_{2}$ to a martingale
difference we shall denote by $D_{0}.$ Moreover $\sup_{1\leq l\leq
m}||\mathbb{E}(S_{l}|\mathcal{F}_{0})||^{2}/m\rightarrow0$. Denote $D_{i}$ the
limit of $D_{i}^{n}$ and construct the martingale $M_{n}=\sum\nolimits_{j=0}%
^{n-1}D_{j}.$

Let $n$ and $m$ be two strictly positive integers. By the fact that both
$D_{0}^{n}$ and $D_{0}^{m}$ are martingale differences and using (\ref{rest})
and (\ref{deco}) we deduce
\[
||D_{0}^{n}-D_{0}^{m}||^{2}=\frac{||M_{m}^{n}-M_{m}^{m}||^{2}}{m}\leq\frac
{1}{m}||(\theta_{0}^{n}-\theta_{m}^{n}+\overline{R}_{m}^{\,n})-(\theta_{0}%
^{m}-\theta_{m}^{m}+\overline{R}_{m}^{m})||^{2}\text{ .}%
\]
So for $n$ fixed, by the fact that $\sup_{1\leq l\leq m}||\mathbb{E}%
(S_{l}|\mathcal{F}_{0})||^{2}/m\rightarrow0$ we have that
\begin{equation}
||D_{0}^{n}-D_{0}||=\lim_{m\rightarrow\infty}||D_{0}^{n}-D_{0}^{m}||\leq
\lim_{m\rightarrow\infty}\frac{1}{m^{1/2}}||\overline{R}_{m}^{n}||=||Y_{0}%
^{n}||_{+}\text{ .} \label{limD}%
\end{equation}
We continue the estimate in the following way%
\begin{align*}
\frac{||S_{n}-M_{n}||^{2}}{n}  &  \leq2(\frac{||S_{n}-M_{n}^{n}||^{2}}%
{n}+\frac{||M_{n}^{n}-M_{n}||^{2}}{n})\\
&  \leq2(\frac{||S_{n}-M_{n}^{n}||^{2}}{n}+||D_{0}^{n}-D_{0}||^{2})\text{ .}%
\end{align*}
The lemma follows by combining the estimates in (\ref{estRn}) and
(\ref{limD}). $\lozenge$

\bigskip

Next we estimate $||Y_{0}^{n}||_{+}.$

\begin{lemma}
\label{plus}Under the conditions of Proposition \ref{approx}, for every
$n\geq1$ and any $m\geq1$, we have
\begin{equation}
\frac{1}{n^{1/2}}||\max_{1\leq j\leq n}|\sum_{k=0}^{j-1}Y_{k}^{m}|\text{
}||\ll\sum_{k=m+1}^{\infty}\frac{||\mathbb{E}_{0}(S_{k})||}{k^{3/2}}\text{ }
\label{norm}%
\end{equation}
and
\begin{equation}
||Y_{0}^{m}||_{+}\ll\sum_{k\geq m}\frac{||\mathbb{E}_{0}(S_{k})||}{k^{3/2}%
}\text{ .} \label{plusnorm}%
\end{equation}

\end{lemma}

\textbf{Proof.} In order to prove this inequality, we apply the maximal
inequality in Peligrad and Utev (2005) to the stationary sequence $Y_{0}^{m}$
defined by (\ref{defY}), where $m\leq n$. Then,
\[
||\max_{1\leq j\leq n}|\sum_{k=0}^{j-1}Y_{k}^{m}|\text{ }||\ll n^{1/2}%
(||Y_{0}^{m}||+\Delta(Y_{0}^{m}))\text{ .}%
\]
We estimate now $\Delta(Y_{0}^{m}).$ It is convenient to use the decomposition%
\begin{gather*}
\Delta(Y_{0}^{m})=\sum_{k=1}^{\infty}\frac{1}{k^{3/2}}||\mathbb{E}_{0}%
(Y_{0}^{m}+...+Y_{k-1}^{m})||\leq\\
\sum_{k=1}^{m}\frac{1}{k^{3/2}}||\mathbb{E}_{0}(Y_{0}^{m}+...+Y_{k-1}%
^{m})||+\sum_{k=m+1}^{\infty}\frac{1}{k^{3/2}}||\mathbb{E}_{0}(Y_{0}%
^{m}+...+Y_{k-1}^{m})||\text{ .}%
\end{gather*}
To estimate the first sum notice that, by the properties of conditional
expectation, we have
\[
||\mathbb{E}_{0}(Y_{0}^{m}+...+Y_{k-1}^{m})||\leq k||\mathbb{E}_{0}(Y_{0}%
^{m})||\,,
\]
and then, since $||\mathbb{E}_{0}(Y_{0}^{m})||\leq||\mathbb{E}_{0}(S_{m})||/m$
we have%
\[
\sum_{k=1}^{m}\frac{1}{k^{3/2}}||\mathbb{E}_{0}(Y_{0}^{m}+...+Y_{k-1}%
^{m})||\leq\frac{1}{m}\sum_{k=1}^{m}\frac{||\mathbb{E}_{0}(S_{m})||}{k^{1/2}%
}\ll\frac{1}{m^{1/2}}||\mathbb{E}_{0}(S_{m})||\text{ .}%
\]
To estimate the second sum we also apply the properties of conditional
expectation and write this time%
\[
||\mathbb{E}_{0}(Y_{0}^{m}+...+Y_{k-1}^{m})||\leq||\mathbb{E}_{0}%
(S_{k})||\text{ .}%
\]
Then,%
\[
\sum_{k=m+1}^{\infty}\frac{1}{k^{3/2}}||\mathbb{E}_{0}(Y_{0}^{m}%
+...+Y_{k-1}^{m})||\leq\sum_{k=m+1}^{\infty}\frac{||\mathbb{E}_{0}(S_{k}%
)||}{k^{3/2}}\,,
\]
and overall, for a certain positive constant $C,$%
\[
\Delta(Y_{0}^{m})\leq C{\Large (}\frac{1}{m^{1/2}}||\mathbb{E}_{0}%
(S_{m})||+\sum_{k=m+1}^{\infty}\frac{||\mathbb{E}_{0}(S_{k})||}{k^{3/2}%
}{\Large )}\text{ .}%
\]
We conclude that for any strictly positive integers $n$ and $m\ $ we have%
\[
\frac{1}{\sqrt{n}}{\large ||}\max_{1\leq j\leq n}|\sum_{k=0}^{j-1}Y_{k}%
^{m}|{\large ||}_{2}\leq2\frac{||\mathbb{E}_{0}(S_{m})||}{m}+80C{\Large (}%
\frac{||\mathbb{E}_{0}(S_{m})||}{m^{1/2}}+\sum_{k=m+1}^{\infty}\frac
{||\mathbb{E}_{0}(S_{k})||}{k^{3/2}}{\Large )}\text{ .}%
\]
$\ $\ The estimate (\ref{norm}) of this lemma follows now by using Lemma
\ref{SUB} from the Appendix with $p=2$ and $\gamma=1/2$. With the notation
(\ref{norm+}), by passing to the limit in relation (\ref{norm}), we obtain
(\ref{plusnorm}). $\ \lozenge$

\bigskip

\textbf{Proof of Proposition \ref{approx}}.

\bigskip

Notice that (\ref{MW}) implies $||Y_{0}^{m}||_{+}\rightarrow0.$ We combine the
estimate in Lemma \ref{estplus} with the estimate of $\ ||Y_{0}^{m}||_{+}$ in
Lemma \ref{plus} to obtain the desired result, via Lemma \ref{SUB} in Appendix
applied with $p=2$ and $\gamma=1/2$.

\section{Almost sure martingale approximations}

In this section we use the estimate (\ref{Approx}) obtained in Proposition
\ref{approx} to approximate a partial sum by a martingale in the almost sure sense.

\begin{proposition}
\label{almost1}Assume $(b_{n})_{n\geq1}$ is any nondecreasing positive, slowly
varying sequence such that
\begin{equation}
\sum_{n\geq1}\frac{b_{n}}{n}\big (\sum_{k\geq n}\frac{||\mathbb{E}_{0}%
(S_{k})||}{k^{3/2}}\big )^{2}<\infty\text{ .} \label{slow}%
\end{equation}
Then, there is a martingale $(M_{k})_{k\in\mathbb{Z}}$ with stationary and
square integrable differences adapted to $(\mathcal{F}_{k+1})_{k\in\mathbb{Z}%
}$ satisfying
\begin{equation}
\frac{S_{n}-M_{n}}{\sqrt{nb_{n}^{\ast}}}\rightarrow0\text{ a.s.} \label{AS1}%
\end{equation}
where $b_{n}^{\ast}:=\sum\nolimits_{k=1}^{n}(kb_{k})^{-1}$.
\end{proposition}

As an immediate consequence of this proposition we formulate the following corollary:

\begin{corollary}
\label{QUENCED}Assume that for a certain sequence of positive numbers
$(b_{n})_{n\geq1}$ that is slowly varying, nondecreasing and satisfies
$\sum\nolimits_{n\geq1}(nb_{n})^{-1}<\infty$, condition (\ref{slow}) is
satisfied. Then there is a martingale $(M_{k})_{k\in\mathbb{Z}}$ with
stationary and square integrable differences adapted to $(\mathcal{F}%
_{k+1})_{k\in\mathbb{Z}}$ satisfying:%
\begin{equation}
\frac{S_{n}-M_{n}}{n^{1/2}}\rightarrow0\text{ a.s.} \label{AS}%
\end{equation}

\end{corollary}

\textbf{Example: }In Corollary \ref{QUENCED} the sequence $(b_{n})_{n\geq3}$
can be taken for instance\textbf{ }$\ b_{n}=(\log n)(\log_{2}n)^{\gamma}%
$\textbf{ }for some\textbf{ }$\gamma>1.$

\bigskip

Selecting in Proposition \ref{almost1} the sequence $b_{n}=\log n,$ we obtain:

\begin{corollary}
\label{LOG}Assume that
\begin{equation}
\sum_{n\geq1}\frac{\log n}{n}\big (\sum_{k\geq n}\frac{||\mathbb{E}_{0}%
(S_{k})||}{k^{3/2}}\big )^{2}<\infty\text{ .} \label{MWlog2}%
\end{equation}
Then there is a martingale $(M_{k})_{k\in\mathbb{Z}}$ with stationary and
square integrable differences adapted to $(\mathcal{F}_{k+1})_{k\in\mathbb{Z}%
}$ satisfying:%
\begin{equation}
\frac{S_{n}-M_{n}}{(n\log_{2}n)^{1/2}}\rightarrow0\text{ a.s.} \label{LOGAS}%
\end{equation}

\end{corollary}

\textbf{Proof of Proposition \ref{almost1}. }

\bigskip

By Corollary 4.2 in Cuny (2009 a), given in Appendix for the convenience of
the reader (see Proposition \ref{Cunya.s.}), in order to show that%
\[
\frac{S_{n}-M_{n}}{\sqrt{nb_{n}^{\ast}}}\rightarrow0\text{ a.s. ,}%
\]
we have to verify that
\[
\sum_{n\geq1}\frac{b_{n}||S_{n}-M_{n}||^{2}}{n^{2}}<\infty\text{ .}%
\]
By Proposition \ref{approx} we know that%

\[
\frac{||S_{n}-M_{n}||}{n^{1/2}}\ll\sum_{k\geq n}\frac{||\mathbb{E}_{0}%
(S_{k})||}{k^{3/2}}\text{ .}%
\]
Therefore the condition (\ref{slow}) implies the desired martingale
approximation. $\lozenge$

\bigskip

\textbf{Remark.} Notice that our condition (\ref{slow}) is implied by the
condition in Corollary 5.8 in Cuny (2009 a). He assumed for the same results
$||\mathbb{E}_{0}(S_{n})||\ll n^{1/2}(\log n)^{-2}(\log_{2})^{-1-\delta}$ with
$\delta>0,$ that clearly implies (\ref{slow}). Also (\ref{MWlog2}) is implied
by the result in Corollary 5.7 in Cuny (2009 a) who obtained the same result
under the condition $||\mathbb{E}_{0}(S_{n})||\ll n^{1/2}(\log n)^{-2}%
(\log_{2}n)^{-\tau}$ with $\tau>1/2$.

\bigskip

In the next two subsections we propose two ways to improve on the rate of
convergence to $0$ of $||\mathbb{E}_{0}(S_{k})||/\sqrt{k}$ that assure an
almost sure martingale approximations in some sense.

\subsection{\textbf{Averaging }}

In the next proposition we study a Cezaro-type almost sure martingale approximation.

\begin{proposition}
\label{AVE}Assume that
\begin{equation}
\sum_{n\geq1}\frac{1}{n}\big (\sum_{k\geq n}\frac{||\mathbb{E}_{0}(S_{k}%
)||}{k^{3/2}}\big )^{2}<\infty\text{ .} \label{MWave}%
\end{equation}
Then there is a martingale $(M_{k})_{k\in\mathbb{Z}}$ with stationary and
square integrable differences adapted to $(\mathcal{F}_{k+1})_{k\in\mathbb{Z}%
}$ satisfying:
\begin{equation}
\frac{1}{n}\sum_{k=1}^{n}\frac{|S_{k}-M_{k}|}{k^{1/2}}\rightarrow0\text{ a.s.
} \label{AVE-AS}%
\end{equation}

\end{proposition}

Before proving this proposition we shall formulate condition (\ref{MWave}) in
an equivalent form that is due to monotonicity:%
\begin{equation}
\sum_{r\geq0}\Big (\sum_{l\geq2^{r}}\frac{||\mathbb{E}(S_{l}|\mathcal{F}%
_{0})||}{l^{3/2}}\Big )^{2}<\infty\text{ .} \label{ModMW1}%
\end{equation}

\textbf{Proof of Proposition \ref{AVE}.} We notice that the condition
(\ref{MWave}) implies by Proposition \ref{approx} the existence of a
martingale $(M_{n})_{n\geq0}$ with stationary differences such that
\begin{equation}
\sum_{n\geq1}\frac{||S_{n}-M_{n}||^{2}}{n^{2}}<\infty\text{ ,} \label{cond2}%
\end{equation}
that further implies
\[
\sum_{n\geq1}\frac{(S_{n}-M_{n})^{2}}{n^{2}}<\infty\text{ a.s.}%
\]
Whence, by Kronecker lemma,
\[
\frac{1}{n}\sum_{k=1}^{n}\frac{(S_{k}-M_{k})^{2}}{k}\rightarrow0\text{ a.s. }%
\]
and then, by Cauchy-Schwarz inequality
\[
\sum_{k=1}^{n}\frac{|S_{k}-M_{k}|}{k^{1/2}}\leq\big (n\sum_{k=1}^{n}%
\frac{(S_{k}-M_{k})^{2}}{k} \big )^{1/2}\text{ .}%
\]
Therefore%
\[
\frac{1}{n}\sum_{k=1}^{n}\frac{|S_{k}-M_{k}|}{k^{1/2}}\rightarrow0\text{
}\ \ \text{a.s. }%
\]
$\lozenge$

\bigskip

This idea\textbf{ }of considering the average approximation can be also
applied to Markov chains with normal operators (i.e. $QQ^{\ast}=Q^{\ast}Q$ on
{$\mathbb{L}^{2}(\pi)$)}. For this case we can replace our Proposition
\ref{approx} by a result stated in Cuny (2009, a) for normal Markov chains
namely, { }
\begin{equation}
\frac{||S_{n}-M_{n}||^{2}}{n}\ll\frac{1}{n}\sum_{k\leq n}\frac{||\mathbb{E}%
_{0}(S_{k})||^{2}}{k}+\sum_{k>n}\frac{||\mathbb{E}_{0}(S_{k})||^{2}}{k^{2}%
}\text{ .} \label{Cuny}%
\end{equation}
Then we can replace in the proof of Proposition \ref{AVE} our Proposition
\ref{approx} by the inequality (\ref{Cuny}).

Notice that
\[
\sum_{n\geq1}\frac{1}{n}\big (\frac{1}{n}\sum_{k\leq n}\frac{||\mathbb{E}%
_{0}(S_{k})||^{2}}{k} \big )\ll\sum_{k\geq1}\frac{||\mathbb{E}_{0}%
(S_{k})||^{2}}{k^{2}}%
\]
and%
\[
\sum_{n\geq1}\frac{1}{n}\big (\sum_{k\geq n}\frac{||\mathbb{E}_{0}%
(S_{k})||^{2}}{k^{2}} \big )\ll\sum_{n\geq1}\frac{\log n||\mathbb{E}_{0}%
(S_{n})||^{2}}{n^{2}}<\infty\text{ .}%
\]
We can then formulate:

\begin{proposition}
Given a stationary Markov chain $(\xi_{n})_{n\in\mathbb{Z}}$ with Normal
operator and $X_{i}=f(\xi_{i})$ is centered at mean and square integrable. If
the condition%
\begin{equation}
\sum_{n\geq2}\frac{\log n||\mathbb{E}_{0}(S_{n})||^{2}}{n^{2}}<\infty\, ,
\label{Nornal}%
\end{equation}
is satisfied, then (\ref{AVE-AS}) holds.
\end{proposition}

We point out that condition (\ref{Nornal}) by itself does not imply (\ref{AS})
so the averaging is needed. As a matter of fact, Cuny and Peligrad (2009)
commented that there is a stationary and ergodic normal Markov chain and a
function $f$ such that
\[
\sum_{n\geq2}\frac{\log n\log_{2}n||\mathbb{E}_{0}(S_{n})||^{2}}{n^{2}}%
<\infty\, ,
\]
and such that (\ref{AS}) fails.

\subsection{\textbf{Higher moments}}

Another way to improve on the rate of convergence to $0$ of $||\mathbb{E}%
(S_{j}|\mathcal{F}_{0})||/j^{1/2}$ in order to establish central limit
theorems started at a point is to consider the existence of moments higher
than $2$.

\begin{proposition}
\label{QUENCHED-H}\textit{Assume that }for some $\delta>0$, $\mathbb{E}%
|X_{0}|^{2+\delta}<\infty,$ $\ $and that the condition (\ref{MWave}) is
satisfied. \textit{Then, }there is a martingale $(M_{k})_{k\in\mathbb{Z}}$
with stationary and square integrable differences adapted to $(\mathcal{F}%
_{k+1})_{k\in\mathbb{Z}}$ satisfying \textit{for every }$\varepsilon>0$
\textit{\ }%
\[
\sum_{n\geq1}n^{-1}\mathbb{P}(\max_{j\leq n}|S_{j}-M_{j}|\geq\varepsilon
\sqrt{n})<\infty\,,
\]
and therefore $S_{n}-M_{n}=o(n^{1/2})$ a.s.
\end{proposition}

\bigskip

\textbf{Proof of Proposition \ref{QUENCHED-H}. } Since, by assumption
(\ref{MWave}) it follows that $\sum_{j\geq1}j^{-3/2}||\mathbb{E}%
(S_{j}|\mathcal{F}_{0})||<\infty,$ according to Proposition \ref{approx},
there exists a martingale $(M_{k})_{k\in\mathbb{Z}}$ with stationary and
square integrable differences $(D_{k})_{k\in\mathbb{Z}}$ adapted to
$(\mathcal{F}_{k+1})_{k\in\mathbb{Z}}$ such that (\ref{Approx}) is satisfied.
Applying then Corollary \ref{cormaxstat} with $\varphi(x)=x^{2}$, $p=2$,
$Y_{i}=X_{i-1}$, $Z_{i}=D_{i-1}$ and $\mathcal{G}_{i}=\mathcal{F}_{i}$, we get
that for every $\varepsilon>0$ and any $\alpha\in\lbrack0,1)$,
\begin{align}
\mathbb{P}(\max_{j\leq n}|S_{j}-M_{j}|  &  \geq4\varepsilon\sqrt{n})\ll
\frac{1}{\varepsilon^{2}}\Big (\sum_{j\geq n}\frac{1}{j^{3/2}}||\mathbb{E}%
(S_{j}|\mathcal{F}_{0})||\Big )^{2}+\frac{n^{1/2}}{\varepsilon}{\mathbb{E}%
}\left(  |X_{0}|\mathbf{1}_{\{|X_{0}|\geq\varepsilon n^{1/2-\alpha}\}}\right)
\label{ine1}\\
&  +\frac{1}{\varepsilon^{2}} \Big (\sum_{k\geq\lbrack n^{\alpha}]+1}%
\frac{||\mathbb{E}(S_{l}|\mathcal{F}_{0})||}{k^{3/2}}\Big )^{2}\,.\nonumber
\end{align}
Choosing now $\alpha=\delta/(2+2\delta)$, we get by using Fubini theorem that
\begin{equation}
\sum_{n\geq1}\frac{1}{n^{1/2}}{\mathbb{E}}\left(  |X_{0}|\mathbf{1}%
_{\{|X_{0}|\geq\varepsilon n^{1/2-\alpha}\}}\right)  \ll\frac{1}%
{\varepsilon^{1+\delta}}{\mathbb{E}}|X_{0}|^{2+\delta}\,. \label{conv1}%
\end{equation}
Therefore, starting from (\ref{ine1}) and using (\ref{conv1}), we infer that
the theorem holds provided that
\begin{equation}
\sum_{n\geq1}\frac{1}{n}\Big (\sum_{j\geq\lbrack n^{\delta/(2+2\delta)}]}%
\frac{||\mathbb{E}(S_{j}|\mathcal{F}_{0})||}{j^{3/2}}\Big )^{2}<\infty\,.
\label{conv3}%
\end{equation}
Now by the usual comparison between the series and the integrals, we notice
that for any nonincreasing and positive function $h$ on ${\mathbb{R}}^{+}$ and
any positive $\alpha$, $\sum_{n\geq1}n^{-1}h(n^{\alpha})<\infty$ if and only
if $\sum_{n\geq1}n^{-1}h(n)<\infty$. Applying this result with
$h(x)=\Big (\sum_{j\geq\lbrack x]}\frac{||\mathbb{E}(S_{j}|\mathcal{F}_{0}%
)||}{j^{3/2}}\Big )^{2}$, it follows that the conditions (\ref{MWave}) and
(\ref{conv3}) are equivalent. This ends the proof of the theorem. $\lozenge$

\bigskip

Next proposition will be useful to transport from the martingale to the
stationary sequence the law of iterated logarithm. Its proof follows the same
line as of Proposition \ref{QUENCHED-H} with the obvious changes.

\begin{proposition}
\label{LOG-H}\textit{Assume that }$\mathbb{E}|X_{0}|^{2+\delta}<\infty$ for
some $\delta>0$, and that
\begin{equation}
\sum_{n\geq3}\frac{1}{n\log_{2}n}{\large (}\sum_{j\geq n}\frac{||\mathbb{E}%
(S_{j}|\mathcal{F}_{0})||}{j^{3/2}}{\large )}^{2}<\infty\text{ .}
\label{highlog}%
\end{equation}
Then there is a martingale $(M_{k})_{k\in\mathbb{Z}}$ with stationary and
square integrable differences adapted to $(\mathcal{F}_{k+1})_{k\in\mathbb{Z}%
}$ such that we have the approximation $S_{n}-M_{n}=o(n\log_{2}n)^{1/2}$ a.s..
\end{proposition}

\bigskip

We shall point now two sets of conditions that satisfy the conditions of these
last two propositions. Assume that $||\mathbb{E}(S_{n}|\mathcal{F}_{0})||\ll
n^{1/2}(\log n)^{-3/2}(\log_{2}n)^{\beta}$ for a certain $\beta>1/2.$ Then
condition (\ref{MWave}) is satisfied. If $||\mathbb{E}(S_{n}|\mathcal{F}%
_{0})||\ll n^{1/2}(\log n)^{-3/2}(\log_{2}n)^{-\gamma}$ for a $\gamma>0$ then
the condition (\ref{highlog}) is satisfied.

\section{Quenched Central Limit Theorem and the Law of iterated Logarithm}

We shall formulate here a few applications of the almost sure martingale
approximations to quenched functional CLT and LIL. For simplicity we assume in
this section that the stationary sequence is ergodic to avoid random normalizers.

\begin{theorem}
Assume that the stationary sequence is ergodic and the conditions of Corollary
\ref{QUENCED} or Proposition \ref{QUENCHED-H} hold. Then
\[
S_{[nt]}/\sqrt{n}\Rightarrow W(\sigma^{2}t)\text{ under }\mathbb{P}_{0}\text{
a.s.} \, ,
\]
where $\sigma=||D_{0}||$ and $D_{0}$ is defined by (\ref{martd}).
\end{theorem}

\textbf{Proof}. The conditions of Corollary \ref{QUENCED} or Proposition
\ref{QUENCHED-H} imply that for every $\varepsilon>0$
\[
\mathbb{P}_{0}(\max_{1\leq k\leq n}|S_{k}-M_{k}|>\varepsilon\sqrt
{n})\rightarrow0\text{ \ \ a.s.}%
\]
that further implies
\[
\mathbb{P}_{0}(\sup_{0\leq t\leq1}|S_{[nt]}-M_{[nt]}|>\varepsilon\sqrt
{n})\rightarrow0\text{ }\ \ \text{a.s.}%
\]
According to Theorem 3.1 in Billingsley (1999), the limiting distribution of
$S_{[nt]}|/\sqrt{n}$ is the same as of $M_{[nt]}/\sqrt{n}$ under
$\mathbb{P}_{0}$ a.s. It was shown in Derriennic and Lin (2001, a) in details
that
\[
M_{[nt]}/\sqrt{n}\Rightarrow W(\sigma^{2}t)\text{ under }\mathbb{P}_{0}\text{
a.s.}\,,
\]
and the result follows. $\lozenge$

\begin{theorem}
\label{CLTave}Assume that the stationary sequence is ergodic and the
conditions of Proposition \ref{AVE} are satisfied. Then we have
\begin{equation}
\frac{1}{n}\sum_{k=1}^{n}\frac{S_{k}}{k^{1/2}}\Rightarrow\sigma N(0,\frac
{2}{3})\text{ \ \ under }\mathbb{\ \mathbb{P}}_{0}\text{ a.s. },
\label{CLTave1}%
\end{equation}
where $\sigma=||D_{0}||$ and $D_{0}$ is defined by (\ref{martd}).
\end{theorem}

\textbf{Proof}. Under the condition (\ref{MWave}) we know there is an ergodic
martingale $(M_{n})$ with stationary and square integrable differences
$(D_{n})$ satisfying
\[
\frac{1}{n}\sum_{k=1}^{n}\frac{|S_{k}-M_{k}|}{k^{1/2}}\rightarrow0\text{
a.s.}\,.
\]
Then, by Theorem 3.1 in Billingsley (1999), the limiting distribution of
\[
\frac{1}{n}\sum_{k=1}^{n}\frac{S_{k}}{k^{1/2}}\text{ coincides to the limiting
distribution of }\frac{1}{n}\sum_{k=1}^{n}\frac{M_{k}}{k^{1/2}}\text{ under
}\mathbb{P}_{0}\text{ a.s.}%
\]
By changing the order of summation we can rewrite $\sum_{k=1}^{n}M_{k}%
/k^{1/2}$ as
\[
\sum_{i=0}^{n-1}(\sum_{k=i+1}^{n}\frac{1}{k^{1/2}})D_{i}%
\]
and, according to the Raikov method for proving the central limit theorem for
martingales, we have to study the limit of the sum of squares. Then, starting
from
\[
\frac{1}{n}\sum_{i=0}^{n-1}D_{i}^{2}\rightarrow\mathbb{E}(D_{0}^{2}%
)=\sigma^{2}\text{ \ \ a.s. and in }\mathbb{L}_{1}\text{,}%
\]
we easily conclude, by the generalized Toeplitz lemma (see Lemma \ref{GKL})
that
\[
\ \frac{1}{n^{2}}\sum_{i=0}^{n-1}(\sum_{k=i+1}^{n}\frac{1}{k^{1/2}})^{2}%
D_{i}^{2}\rightarrow\frac{2}{3}\sigma^{2}\text{ \ \ a.s. and in }%
\mathbb{L}_{1}\text{.}%
\]
Then, by Theorem 3.6 in Hall and Heyde (1980) we easily obtain the convergence
in (\ref{CLTave1}). $\ $ $\lozenge$

\begin{theorem}
Assume that either the conditions of Corollary \ref{LOG} or of Proposition
\ref{LOG-H} hold and in addition the sequence is ergodic. Then we can redefine
$(X_{n})_{n\in\mathbb{Z}},$ without changing its distribution, on a richer
probability space on which there exists a standard Brownian motion
$(W(t),t\geq0)$ such that\textit{ }%
\[
S_{n}-W(n||D_{0}||^{2})=o(n\log_{2}n)^{1/2}\text{ a.s.}%
\]
Therefore, the LIL holds:\
\[
\lim\sup\pm\frac{S_{n}}{(2n\log_{2}n)^{1/2}}=||D_{0}||\text{ a.s.}%
\]

\end{theorem}

\textbf{Proof}. Since by Corollary \ref{LOG} $\ $or by Proposition \ref{LOG-H}
$\ $we have $S_{n}-M_{n}=o(n\log_{2}n)^{1/2}$ $a.s.$ the result follows by the
almost sure invariance principle for stationary and square integrable
martingales (see Strassen, 1967). $\lozenge$

\section{Applications}

We shall mention two examples for which the quantity $||\mathbb{E}%
(S_{j}|\mathcal{F}_{0})||^{2}$ is estimated. Then, these estimates introduced
in our results will provide new asymptotic results started at a point and LIL,
that improve the previous results in the literature.

\bigskip

\textbf{1.} \textbf{Application to linear processes.}

\bigskip

Let $(\varepsilon_{n})_{n\in\mathbb{Z}}$ be a sequence of ergodic martingale
differences and consider the linear process
\[
X_{k}=\sum_{i\geq1}a_{i}\varepsilon_{k-i}%
\]
where $(a_{i})_{i\geq1}$ is a sequence of real constants such that
$\sum_{i\geq1}a_{i}^{2}<\infty$. We define
\[
S_{n}=\sum_{i=1}^{n}X_{i}\text{ .}%
\]
Denote by
\[
b_{nj}=a_{j+1}+...+a_{j+n}\text{ .}%
\]
Then
\[
||\mathbb{E}(S_{n}|\mathcal{F}_{0})||^{2}=\sum_{j\geq0}b_{nj}^{2}\text{ .}%
\]
For the particular case $a_{n}\ll1/(nL(n))$, where $L($\textperiodcentered$)$
is a positive, nondecreasing, slowly varying function, computations in Zhao
and Woodroofe (2008-a) show that $||\mathbb{E}(S_{n}|\mathcal{F}_{0}%
)||\ll\sqrt{n}/L(n)$.

\bigskip

\textbf{2}. \textbf{Application to functions of Bernoulli shifts.}

\bigskip

Let $({\varepsilon}_{k})_{k\in\mathbb{Z}}$ be an i.i.d. sequence of Bernoulli
variables, that is $\mathbb{P}({\varepsilon}_{1}=0)=1/2=\mathbb{P}%
({\varepsilon}_{1}=1)$ and let
\[
Y_{n}=\sum_{k=0}^{\infty}2^{-k-1}{\varepsilon}_{n-k},\quad X_{n}=g(Y_{n}%
)-\int_{0}^{1}g(x)dx,\text{\ and }S_{n}=\sum_{k=1}^{n}X_{k}\text{ ,}%
\]
where $g\in\mathbb{L}_{2}(0,1)$, $(0,1)$ being equipped with the Lebesgue
measure. The transform $Y_{j}$ is usually referred to as the Bernoulli shift
of the i.i.d. sequence $({\varepsilon}_{k})_{k\in\mathbb{Z}}$. Then, following
Maxwell and Woodroofe (2000), as in Peligrad et al (2007),%

\[
||{\mathbb{E}}(g(Y_{k})|Y_{0})||^{2}\leq2^{k}\int_{0}^{1}\int_{0}%
^{1}I_{(|x-y|\leq2^{-k})}|g(x)-g(y)|^{2}dydx \, ,
\]
and then
\[
||\mathbb{E}(S_{j}|\mathcal{F}_{0})||^{2}\leq\sum_{k=1}^{j}2^{k}\int_{0}%
^{1}\int_{0}^{1}I_{(|x-y|\leq2^{-k})}|g(x)-g(y)|^{2}dydx\text{ .}%
\]

\section{Appendix}

\subsection{Maximal inequalities}

Following the idea of proof of the maximal inequality given in Proposition 5
of Merlev\`{e}de and Peligrad (2010), we shall prove the following result:

\begin{proposition}
\label{propmaxineproba} Let $(Y_{i})_{1\leq i\leq2^{r}}$ be real random
variables where $r$ is a positive integer. Assume that the random variables
are adapted to an increasing filtration $({\mathcal{G}}_{i})_{1\leq i\leq
2^{r}}$. Let $(Z_{i})_{1\leq i\leq2^{r}}$ be real random variables adapted to
$({\mathcal{G}}_{i})_{1\leq i\leq2^{r}}$ and such that for every $i$,
${\mathbb{E}}(Z_{i}|{\mathcal{G}}_{i-1})=0$ a.s. Let $S_{n}=Y_{1}+\cdots
+Y_{n}$ and $T_{n}=Z_{1}+\cdots+Z_{n}$. Let $\varphi$ be a nondecreasing, non
negative, convex and even function. Then for any positive real $x$, any real
$p\geq1$ and any integer $u\in\lbrack0,r-1]$, the following inequality holds:
\begin{align*}
&  {\mathbb{P}}(\max_{1\leq i\leq2^{r}}|S_{i}-T_{i}|\geq4x)\leq\frac
{1}{\varphi(x)}{\mathbb{E}}(\varphi(S_{2^{r}}-T_{2^{r}}))+\frac{1}{x}%
\sum_{i=1}^{2^{r}}{\mathbb{E}}(|Y_{i}|\mathbf{1}_{\{|Y_{i}|\geq x/2^{u}\}})\\
&  \quad\quad\quad\quad+\frac{1}{x^{p}}\ \Big(\sum_{l=u}^{r-1}\Big (\sum
_{k=1}^{2^{r-l}-1}||{\mathbb{E}}(S_{(k+1)2^{l}}-S_{k2^{l}}|{{\mathcal{G}}%
}_{k2^{l}})||_{p}^{p}\Big )^{1/p}\Big )^{p}\,.
\end{align*}

\end{proposition}

\begin{remark}
When the sequence $(Y_{n})_{n\in\mathbb{Z}}$ is stationary as well as the
filtration $({\mathcal{G}}_{n})_{n\in\mathbb{Z}}$, the inequality has the
following form:
\begin{align*}
&  {\mathbb{P}}(\max_{1\leq i\leq2^{r}}|S_{i}-T_{i}|\geq4x)\leq\frac
{1}{\varphi(x)}{\mathbb{E}}(\varphi(S_{2^{r}}-T_{2^{r}}))+\frac{2^{r}}%
{x}{\mathbb{E}}(|Y_{1}|\mathbf{1}_{\{|Y_{1}|\geq x/2^{u}\}})\\
&  \quad\quad\quad\quad+\frac{1}{x^{p}}2^{r}\Big(\sum_{l=u}^{r-1}\frac
{1}{2^{l/p}}||{\mathbb{E}}(S_{2^{l}}|{{\mathcal{G}}}_{0})||_{p}\Big )^{p}\,.
\end{align*}

\end{remark}

Notice now that for any integer $n\in\lbrack2^{r-1},2^{r})$, where $r$ is a
positive integer, ${\mathbb{E}}(\varphi(S_{2^{r}}-T_{2^{r}}))\leq\max
_{n<k<2n}{\mathbb{E}}(\varphi(S_{k}-T_{k}))$. In addition, due to the
subadditivity of the sequence $\big (||{\mathbb{E}}(S_{n}|{\mathcal{G}}%
_{0})||_{p}\big)_{n\geq1}$, we have that
\[
2^{k}||{\mathbb{E}}(S_{2^{k}}|{\mathcal{G}}_{0})||_{p}\leq2\sum_{j=1}%
^{k}||{\mathbb{E}}(S_{j}|{\mathcal{G}}_{0})||_{p}\,,
\]
implying that for any integer $n\in\lbrack2^{r-1},2^{r})$, where $r$ is a
positive integer, and any integer $u\in\lbrack0,r-1]$,
\[
\sum_{i=u}^{r-1}\frac{1}{2^{i/p}}||{\mathbb{E}}(S_{2^{i}}|{\mathcal{G}}%
_{0})||_{p}\leq2\sum_{j=1}^{2^{r-1}}||{\mathbb{E}}(S_{j}|{\mathcal{G}}%
_{0})||_{p}\sum_{i:2^{i}\geq j\vee2^{u}}\frac{1}{2^{i(1+1/p)}}\,,
\]
and then that
\[
\sum_{i=u}^{r-1}\frac{1}{2^{i/p}}||{\mathbb{E}}(S_{2^{i}}|{\mathcal{G}}%
_{0})||_{p}\leq\frac{2^{2+1/p}}{2^{1+1/p}-1}\Big (\frac{1}{2^{u(1+1/p)}}%
\sum_{k=1}^{2^{u}-1}||{\mathbb{E}}(S_{k}|{\mathcal{G}}_{0})||_{p}%
+\sum_{k=2^{u}}^{n}\frac{||{\mathbb{E}}(S_{k}|{\mathcal{G}}_{0})||_{p}%
}{k^{1+1/p}}\Big )\,.
\]
It remains to apply Lemma \ref{SUB} below to obtain the following corollary:

\begin{corollary}
\label{cormaxstat} Let $(Y_{i})_{i\in\mathbb{Z}}$ be a stationary sequence of
real random variables. Assume that the random variables are adapted to an
increasing and stationary filtration $({\mathcal{G}}_{i})_{i\in\mathbb{Z}}$.
Let $(Z_{i})_{i\in\mathbb{Z}}$ be a sequence of real random variables adapted
to $({\mathcal{G}}_{i})_{i\in\mathbb{Z}}$ and such that for all $i$,
${\mathbb{E}}(Z_{i}|{\mathcal{G}}_{i-1})=0$ a.s. Let $S_{n}=Y_{1}+\cdots
+Y_{n}$ and $T_{n}=Z_{1}+\cdots+Z_{n}$. Let $\varphi$ be a nondecreasing, non
negative, convex and even function. Then for any positive real $x$, any
positive integer $n$, any real $p\geq1$ and any real $\alpha\in\lbrack0,1]$,
the following inequality holds:
\begin{gather*}
{\mathbb{P}}(\max_{1\leq i\leq n}|S_{i}-T_{i}|\geq4x)\leq\frac{1}{\varphi
(x)}\max_{n<k<2n}{\mathbb{E}}(\varphi(S_{k}-T_{k}))+\frac{2n}{x}{\mathbb{E}%
}(|Y_{1}|\mathbf{1}_{\{|Y_{1}|\geq x/n^{\alpha}\}})\\
+\frac{c_{p}n}{x^{p}}\Big(\sum_{k=[n^{\alpha}]+1}^{\infty}\frac{||{\mathbb{E}%
}(S_{k}|{\mathcal{G}}_{0})||_{p}}{k^{1+1/p}}\Big )^{p}\,,
\end{gather*}
where $c_{p}$ is a positive constant depending only on $p$.
\end{corollary}

\textbf{Proof of Proposition \ref{propmaxineproba}}\textit{.}

\bigskip

Using the fact that ${\mathbb{E}}(T_{n}-T_{k}|{\mathcal{G}}_{k})=0$ for $0\leq
k\leq n$, we get, for any $m\in\lbrack0,2^{r}-1]$, that
\[
S_{2^{r}-m}-T_{2^{r}-m}={\mathbb{E}}(S_{2^{r}}-T_{2^{r}}|{\mathcal{G}}%
_{2^{r}-m})-{\mathbb{E}}(S_{2^{r}}-S_{2^{r}-m}|{\mathcal{G}}_{2^{r}-m})\,.
\]
So,%
\begin{equation}
\max_{1\leq i\leq2^{r}}|S_{i}-T_{i}|\leq\max_{0\leq m\leq2^{r}-1}|{\mathbb{E}%
}(S_{2^{r}}-T_{2^{r}}|{\mathcal{G}}_{2^{r}-m})|+\max_{0\leq m\leq2^{r}%
-1}|{\mathbb{E}}(S_{2^{r}}-S_{2^{r}-m}|{\mathcal{G}}_{2^{r}-m}%
)|\,.\label{dec1max}%
\end{equation}
Since $(\mathbb{E}(S_{2^{r}}-T_{2^{r}}|{\mathcal{G}}_{u})_{u\geq1}$ is a
martingale, we shall use Doob's maximal inequality (see Theorem 2.1 in Hall
and Heyde, 1980) to deal with the first term in the right hand side of
(\ref{dec1max}). Hence, since $\varphi$ is a nondecreasing, non negative,
convex and even function, we get that
\begin{equation}
{\mathbb{P}}\big (\max_{0\leq m\leq2^{r}-1}|{\mathbb{E}}(S_{2^{r}}-T_{2^{r}%
}|{\mathcal{G}}_{2^{r}-m})|\geq x\big )\leq\frac{1}{\varphi(x)}{\mathbb{E}%
}(\varphi(S_{2^{r}}-T_{2^{r}}))\,.\label{boundphi}%
\end{equation}
Write now $m$ in basis $2$ as follows:
\[
m=\sum_{i=0}^{r-1}b_{i}(m)2^{i},\ \text{ where $b_{i}(m)=0$ or $b_{i}(m)=1$%
}\,.
\]
Set $m_{l}=\sum_{i=l}^{r-1}b_{i}(m)2^{i}$. With this notation $m_{0}=m$. Let
$0\leq u\leq r-1$ and write that
\[
|{\mathbb{E}}(S_{2^{r}}-S_{2^{r}-m}|{\mathcal{G}}_{2^{r}-m})|\leq|{\mathbb{E}%
}(S_{2^{r}-m_{u}}-S_{2^{r}-m}|{\mathcal{G}}_{2^{r}-m})|+|{\mathbb{E}}%
(S_{2^{r}}-S_{2^{r}-m_{u}}|{\mathcal{G}}_{2^{r}-m})|\,.
\]
Notice first that
\begin{align*}
&  {\mathbb{P}}\big (\max_{0\leq m\leq2^{r}-1}|{\mathbb{E}}(S_{2^{r}-m_{u}%
}-S_{2^{r}-m}|{\mathcal{G}}_{2^{r}-m})|\geq2x\big )\\
&  \quad\quad\quad\quad\leq{\mathbb{P}}\big (\max_{0\leq m\leq2^{r}-1}%
\sum_{j=2^{r}-m+1}^{2^{r}-m_{u}}|{\mathbb{E}}(Y_{j}|{\mathcal{G}}_{2^{r}%
-m})|\geq2x\big )\,.
\end{align*}
Therefore since $|m-m_{u}|\leq2^{u}$,
\begin{align*}
&  {\mathbb{P}}\big (\max_{0\leq m\leq2^{r}-1}|{\mathbb{E}}(S_{2^{r}-m_{u}%
}-S_{2^{r}-m}|{\mathcal{G}}_{2^{r}-m})|\geq2x\big )\\
&  \quad\quad\quad\quad\leq{\mathbb{P}}\big (\max_{0\leq m\leq2^{r}-1}%
\sum_{j=2^{r}-m+1}^{2^{r}-m_{u}}{\mathbb{E}}(|Y_{j}|\mathbf{1}_{\{|Y_{j}|\geq
x/2^{u}\}}|{\mathcal{G}}_{2^{r}-m})\geq x\big )\\
&  \quad\quad\quad\quad\leq{\mathbb{P}}\big (\max_{0\leq m\leq2^{r}-1}%
\sum_{j=1}^{2^{r}}{\mathbb{E}}(|Y_{j}|\mathbf{1}_{\{|Y_{j}|\geq x/2^{u}%
\}}|{\mathcal{G}}_{2^{r}-m})\geq x\big )\,.
\end{align*}
Noticing then that $\big (\sum_{j=1}^{2^{r}}{\mathbb{E}}(|Y_{j}|\mathbf{1}%
_{\{|Y_{j}|\geq x/2^{u}\}}|{\mathcal{G}}_{k})\big )_{k\geq1}$ is a martingale,
Doob's maximal inequality implies that
\begin{equation}
{\mathbb{P}}\big (\max_{0\leq m\leq2^{r}-1}|{\mathbb{E}}(S_{2^{r}-m_{u}%
}-S_{2^{r}-m}|{\mathcal{G}}_{2^{r}-m})|\geq2x\big )\leq x^{-1}\sum
_{i=1}^{2^{r}}{\mathbb{E}}(|Y_{i}|\mathbf{1}_{\{|Y_{i}|\geq x/2^{u}%
\}})\,,\label{b2max}%
\end{equation}
On the other hand, following the proof of the proposition 5 in Merlev\`{e}de
and Peligrad (2010), for any $p\geq1$, we get that
\begin{equation}
\max_{0\leq m\leq2^{r}-1}|{\mathbb{E}}(S_{2^{r}}-S_{2^{r}-m_{u}}|{\mathcal{G}%
}_{2^{r}-m})|^{p}\leq\sum_{l=u}^{r-1}\lambda_{l}^{1-p}\max_{0\leq m\leq
2^{r}-1}|{\mathbb{E}}(A_{r,l}|{\mathcal{G}}_{2^{r}-m})|^{p}\,,\label{ine1MP10}%
\end{equation}
where
\[
\lambda_{l}=\frac{\alpha_{l}}{\sum_{l=u}^{r-1}\alpha_{l}}\text{ with }%
\alpha_{l}=\Big (\sum_{k=1}^{2^{r-l}-1}||{\mathbb{E}}(S_{(k+1)2^{l}}%
-S_{k2^{l}}|{\mathcal{G}}_{k2^{l}})||_{p}^{p}\Big )^{1/p}\,,
\]
and
\[
A_{r,l}=\max_{1\leq k\leq2^{r-l},k\text{ odd}}|{\mathbf{E}}(S_{2^{r}%
-(k-1)2^{l}}-S_{2^{r}-k2^{l}}|{\mathcal{G}}_{2^{r}-k2^{l}})|\,.
\]
Notice now that by Jensen's inequality, $|{\mathbb{E}}(A_{r,l}|{\mathcal{G}%
}_{2^{r}-m})|^{p}\leq{\mathbb{E}}(A_{r,l}^{p}|{\mathcal{G}}_{2^{r}-m})$. Hence
starting from (\ref{ine1MP10}), we get that for any $p\geq1$,
\begin{align*}
&  {\mathbb{P}}\big (\max_{0\leq m\leq2^{r}-1}|{\mathbb{E}}(S_{2^{r}}%
-S_{2^{r}-m_{u}}|{\mathcal{G}}_{2^{r}-m})|\geq x\big )\\
&  \quad\quad\quad\quad\leq{\mathbb{P}}\big (\sum_{l=u}^{r-1}\lambda_{l}%
^{1-p}\max_{0\leq m\leq2^{r}-1}{\mathbb{E}}(A_{r,l}^{p}|{\mathcal{G}}%
_{2^{r}-m})\geq x^{p}\big )\,.
\end{align*}
Next, since $({\mathbf{E}}(A_{r,l}^{p}|{\mathcal{G}}_{k}))_{k\geq1}$ is a
martingale, Doob's maximal inequality entails that
\[
{\mathbb{P}}\big (\max_{0\leq m\leq2^{r}-1}|{\mathbb{E}}(S_{2^{r}}%
-S_{2^{r}-m_{u}}|{\mathcal{G}}_{2^{r}-m})|\geq x\big )\leq x^{-p}\sum
_{l=u}^{r-1}\lambda_{l}^{1-p}{\mathbb{E}}(A_{r,l}^{p})\,.
\]
Taking into account the fact that that ${\mathbb{E}}(A_{r,l}^{p})\leq
\alpha_{l}^{p}$ together with the definition of $\alpha_{l}$ and $\lambda_{l}%
$, we then derive that for any $p\geq1$,
\begin{align}
&  {\mathbb{P}}\big (\max_{0\leq m\leq2^{r}-1}|{\mathbb{E}}(S_{2^{r}}%
-S_{2^{r}-m_{u}}|{\mathcal{G}}_{2^{r}-m})|\geq x\big )\nonumber\\
&  \quad\quad\quad\quad\leq x^{-p}\ \Big(\sum_{l=u}^{r-1}\Big (\sum
_{k=1}^{2^{r-l}-1}||{\mathbb{E}}(S_{(k+1)2^{l}}-S_{k2^{l}}|{{\mathcal{G}}%
}_{k2^{l}})||_{p}^{p}\Big )^{1/p}\Big )^{p}\,.\label{b3max}%
\end{align}
Starting from (\ref{dec1max}) and considering the bounds (\ref{boundphi}),
(\ref{b2max}) and (\ref{b3max}), the proposition follows. $\lozenge$

\subsection{Technical results}

\begin{lemma}
\label{SUB}In the context of stationary sequences, for every $\mathbb{\gamma
}>0$, $n\geq1$ and $p\geq1$,
\[
\frac{1}{n^{\mathbb{\gamma}}}\max_{1\leq k\leq n}||\mathbb{E}_{0}(S_{k}%
)||_{p}\leq c_{\gamma}\sum_{k\geq n+1}\frac{1}{k^{\mathbb{\gamma}+1}%
}||\mathbb{E}_{0}(S_{k})||_{p}\text{ ,}%
\]
where $c_{\gamma}=3\times2^{2\gamma+1}(2^{\mathbb{\gamma}+1}+1).$
\end{lemma}

\textbf{Proof}. We consider $k\geq n$ and start from the inequality%
\[
|\mathbb{E}_{0}(S_{n})|\leq|\mathbb{E}_{0}(S_{k+n})|+|\mathbb{E}_{0}%
(S_{k+n}-S_{n})|\text{ .}%
\]
Then, by the properties of conditional expectation and stationarity%
\[
||\mathbb{E}_{0}(S_{n})||_{p}\leq||\mathbb{E}_{0}(S_{k+n})||_{p}%
+||\mathbb{E}_{0}(S_{k})||_{p}\text{ .}%
\]
So, for any $n\geq1$
\begin{gather*}
\frac{1}{n^{\mathbb{\gamma}}}||\mathbb{E}_{0}(S_{n})||_{p}=\frac
{1}{n^{\mathbb{\gamma}+1}}||\mathbb{E}_{0}(S_{n})||_{p}(\sum_{k=n+1}%
^{2n}1)\leq2^{\gamma+1}||\mathbb{E}_{0}(S_{n})||_{p}\sum_{k=n+1}^{2n}\frac
{1}{k^{\mathbb{\gamma}+1}}\\
\leq2^{\gamma+1}\sum_{k=n+1}^{2n}\frac{1}{k^{\mathbb{\gamma}+1}}%
||\mathbb{E}_{0}(S_{k+n})||_{p}+2^{\gamma+1}\sum_{k=n+1}^{2n}\frac
{1}{k^{\mathbb{\gamma}+1}}||\mathbb{E}_{0}(S_{k})||_{p}\\
\leq2^{2\gamma+2}\sum_{k=n+1}^{2n}\frac{1}{(k+n)^{\mathbb{\gamma}+1}%
}||\mathbb{E}_{0}(S_{k+n})||_{p}+2^{\gamma+1}\sum_{k=n+1}^{2n}\frac
{1}{k^{\mathbb{\gamma}+1}}||\mathbb{E}_{0}(S_{k})||_{p}\text{ .}%
\end{gather*}
We easily derive%
\begin{equation}
\frac{1}{n^{\mathbb{\gamma}}}||\mathbb{E}_{0}(S_{n})||_{p}\leq2^{\gamma
+1}(2^{\mathbb{\gamma}+1}+1)\sum_{l=n+1}^{\infty}\frac{1}{l^{^{\mathbb{\gamma
}+1}}}||\mathbb{E}_{0}(S_{l})||_{p}\text{ .} \label{subad}%
\end{equation}
By writing now%
\[
|\mathbb{E}_{0}(S_{k})|\leq|\mathbb{E}_{0}(S_{k+n})|+|\mathbb{E}_{0}%
(S_{k+n}-S_{k})|\text{ ,}%
\]
by stationarity we obtain%
\[
\max_{1\leq k\leq n}||\mathbb{E}_{0}(S_{k})||_{p}\leq\max_{n\leq k\leq
2n}||\mathbb{E}_{0}(S_{k})||_{p}+\ ||\mathbb{E}_{0}(S_{n})||_{p}\leq
3\max_{n\leq k\leq2n}||\mathbb{E}_{0}(S_{k})||_{p}\text{ ,}%
\]
and the result follows by the inequality (\ref{subad}) applied for each $k$,
$n\leq k\leq2n$. $\lozenge$

\bigskip

The following lemma is concerned with a property of subadditive functions.

\begin{lemma}
\textit{Assume that }$\Sigma_{n=1}^{\infty}||\mathbb{E}_{0}(S_{n}%
)||n^{-3/2}<\infty$\textit{. Then, for a universal constant }$c$\textit{ and
integer }$a$\textit{ }
\[
{\frac{1}{\sqrt{a}}}\sum_{k=1}^{\infty}{\frac{||\mathbb{E}_{0}(S_{ak}%
)||}{k^{3/2}}}\leq c \, (\sum_{l\geq a} \frac{||\mathbb{E}_{0}(S_{l}%
)||}{l^{3/2}})\text{ .}%
\]

\end{lemma}

\textbf{Proof}. The proof of this Lemma is an easy consequence of Lemmas
3.1-3.3 and the proof of Lemma 3.4 in Peligrad and Utev (2006) and the fact
that $||\mathbb{E}_{0}(S_{n})||$ is a subadditive sequence. $\lozenge$

\bigskip

Next result we formulate is Corollary 4.2 in Cuny (2009 a).

\begin{proposition}
\label{Cunya.s.}Assume $(X_{n})_{n\in Z}$ is a stationary sequence of square
integrable random variables and $(b_{n})_{n\geq1}$ a positive nondecreasing
slowly varying sequence. Assume%
\[
\sum_{n\geq1}\frac{b_{n}||S_{n}||^{2}}{n^{2}}<\infty\text{ .}%
\]
Then
\[
\frac{S_{n}}{\sqrt{nb_{n}^{\ast}}}\rightarrow0\text{ a.s.}%
\]
where $b_{n}^{\ast}:=\sum\nolimits_{k=1}^{n}(kb_{k})^{-1}$.
\end{proposition}

\bigskip

We give here a generalized Toeplitz lemma, which is Lemma 5 in M. Peligrad and
C. Peligrad (2011).

\begin{lemma}
\label{GKL} Assume $(x_{i})_{i\geq1}$ and $(c_{i})_{i\geq1}$ are sequences of
real numbers such that
\[
\ \frac{1}{n}\sum_{i=1}^{n}x_{i}\rightarrow L\text{ , }nc_{n}\rightarrow
\infty\text{ and }\frac{c_{1}+...+c_{n}}{nc_{n}}\rightarrow C<1\text{.}%
\]
Then,%
\[
\frac{\sum_{i=1}^{n}c_{i}x_{i}}{\sum_{i=1}^{n}c_{i}}\rightarrow L\text{ .}%
\]

\end{lemma}

\section{Acknowledgement}

Magda Peligrad would like to thank Dalibor Voln\'{y} for suggesting this topic
and for numerous very helpful discussions.

\end{document}